\documentclass[twocolumn]{autart}    % Enable this line and disable the 
\usepackage{graphicx}          % Include this line if your 
\usepackage{tabularx,booktabs}
\usepackage{comment}
\usepackage{float}
\usepackage{color}
\usepackage{amsmath,amssymb}
\usepackage{amsfonts}
\usepackage{textcomp}
\usepackage{psfrag}
\usepackage{multimedia}
\usepackage{fancybox}
\usepackage{mathtools}
\definecolor{wheat}{rgb}{0.96,0.87,0.70}
\definecolor{arash}{rgb}{1,0.8,0.8}
\definecolor{mario}{rgb}{0.8,0.8,1}
\definecolor{seb}{rgb}{0.8,1,0.8}
\definecolor{akhil}{rgb}{0.8,0.5,0.8}

\usepackage{acronym}
\usepackage{amsmath}

\DeclareMathOperator*{\argmin}{arg\,min}
\allowdisplaybreaks
\newcommand{\vect}[1]{\ensuremath{\boldsymbol{\mathrm{#1}}}}

\newtheorem{corollary}{Corollary}
\newtheorem{Lemma}{Lemma}

\newtheorem{Assumption}{Assumption}

\newcommand {\cmatr}[2]{\left\{\begin{array}{#1}#2\end{array}\right.}

\begin{document}

\begin{frontmatter}

\title{Optimality Conditions for Model Predictive Control: Rethinking Predictive Model Design } % Title, preferably not more 
% Optimality Conditions for Model Predictive Control: Rethinking Predictive Model Design or\\

% thanks[footnoteinfo]{nothing}

\author[ntnu]{Akhil S Anand}\ead{akhil.s.anand@ntnu.no},    % Add the 
\author[arash]{Arash Bahari Kordabad}\ead{arashbk@mpi-sws.org}, 
\author[mario]{Mario Zanon}\ead{mario.zanon@imtlucca.it}, 
\author[ntnu]{Sebastien Gros}\ead{sebastien.gros@ntnu.no}               % e- (ead) as shown

\address[ntnu]{Department of Engineering Cybernetics (ITK),Norwegian University of Science and Technology, 7491 Trondheim, Norway. E-mail: {\{akhil.s.anand, sebastien.gros\}@ntnu.no}.}  % Please supply                                              
%\address[sebastien]{Department of Engineering Cybernetics (ITK),Norwegian University of Science and Technology, 7491 Trondheim, Norway. E-mail: {sebastien.gros@ntnu.no}.}             % full addresses
\address[arash]{Max Planck Institute for Software Systems, Kaiserslautern, Germany. E-mail: arashbk@mpi-sws.org.}

\address[mario]{IMT School for Advanced Studies, 55100 Lucca, Italy. E-mail:  {mario.zanon@imtlucca.it}.}             % full addresses
             % full addresses

\begin{keyword}                           % Five to ten keywords,  
MPC; MDP; Optimal Control, System Identification.               % chosen from the IFAC 
\end{keyword}                             % keyword list or with the 
                                          % help of the Automatica 
                                          % keyword wizard

\begin{abstract}                          % Abstract of not more than 200 words.
Optimality is a critical aspect of Model Predictive Control (MPC), especially in economic MPC. However, achieving optimality in MPC presents significant challenges, and may even be impossible, due to inherent inaccuracies in the predictive models. Predictive models often fail to accurately capture the true system dynamics, such as in the presence of stochasticity, leading to suboptimal MPC policies. In this paper, we establish the necessary and sufficient conditions on the underlying prediction model for an MPC scheme to achieve closed-loop optimality. Interestingly, these conditions are counterintuitive to the traditional approach of building predictive models that best fit the data. These conditions present a mathematical foundation for constructing models that are directly linked to the performance of the resulting MPC scheme.
\end{abstract}

\end{frontmatter}

\section{Introduction}
%Optimality in model-based control, MPC
Model Predictive Control (MPC) is a state-of-the-art optimal control approach for complex real-world systems. MPC relies on an underlying predictive model of the true system to optimize the control actions \cite{MPCbook}. This model can take various forms, such as a state-space representation, a data-driven model, or even input-output sequences \cite{verheijen2023handbook}. Additionally, the model may be constructed using first principles (white box), a combination of first principles and data (gray box), or purely data-driven methods (black box) \cite{Mesbah2022}.  Similarly, the model can be estimated through different System Identification (SYSID) approaches. Despite differences in model structure or estimation methods, the performance of an MPC scheme is critically dependent on the accuracy of the predictions generated by its predictive model.  Consequently, model accuracy is critical when the closed-loop performance of the MPC-controlled system is important, which is usually the case in Economic Model Predictive Control (EMPC) \cite{darby2012mpc}. 

% Accurately modeling system dynamics becomes increasingly difficult in case of This challenge is further amplified when dealing with stochastic systems. 
Conventionally, significant emphasis has been placed on improving the prediction accuracy of models in MPC using advanced SYSID techniques and more recent Machine Learning (ML) approaches \cite{hewing2020learning}.  A classic approach is to use ML techniques to improve the prediction accuracy of the model using real-world data, typically via ridge regression \cite{wu2019machine}. This approach of building models that best fit the data assumes that improved model accuracy will lead to better closed-loop performance in the MPC scheme. This assumption, however, is not supported by theory nor often verified in practice \cite{piga2019performance,anand2024data}.

% On the condition
Given the inherent difficulty in accurately modeling a stochastic system, a natural question arises: Can MPC, or any model-based optimization scheme, achieve optimal performance on the real system? Recent works have shown that the answer is affirmative: an MPC scheme can attain optimality despite model inaccuracies by appropriately adjusting the stage cost of the MPC \cite{gros2019data,kordabad2023reinforcement}.  However, these results do not specify how to construct predictive models such that an MPC scheme can achieve optimality. In our previous work \cite{anand2024data}, we introduced a sufficient condition on the predictive model for an MPC scheme to achieve optimality. Building on this foundation, we now introduce the necessary and sufficient conditions for the predictive model such that an MPC scheme achieves optimality. 

In Section \ref{sec:theorem}, we outline these (necessary and sufficient) conditions on the predictive model for an MPC scheme to achieve optimal performance on the true system. Interestingly, these conditions are counterintuitive to the conventional data-fitting approach to model construction, opening a new approach to constructing predictive models for MPC. We base our optimality analysis on the discrete infinite-horizon discounted Markov Decision Process (MDP) framework \cite{sutton2018reinforcement}. The necessary and sufficient conditions proposed in this paper can be readily extended to broader definitions, but we will remain within this framework for the sake of brevity.

The rest of the paper is structured as follows: Section~\ref{sec:background} introduces the necessary background on MDPs and MPC. Section \ref{sec:conditons} provides the necessary and sufficient conditions for an MPC predictive model to achieve optimality, along with the corresponding proofs. In Section \ref{sec:discussions}, we discuss the implications of our findings. Finally, Section \ref{sec:conclusion} presents the conclusions.

\section{Backgroud}\label{sec:background}
This section will briefly introduce the MDP framework for optimal control \cite{sutton2018reinforcement}, which serves as the basis for our optimality analysis, followed by the necessary background on MPC.
\subsection{Markov Decision Processes}
 MDPs assume the stochastic system dynamics of the form:
\begin{align}
\label{eq:system}
\vect{s}_+ \sim \rho \left(\,.\,|\,\vect s, \vect a\,\right),
\end{align}
where $\vect s\in \mathbb S \subseteq \mathbb R^n $ and $\vect a \in \mathbb A \subseteq \mathbb R^m$ are a pair of state and action (also called input in the control literature), respectively, and $\vect{s}_+\in \mathbb S$ is the successor state. The conditional probability density, or more broadly, conditional probability measure, is denoted by $\rho$. For a given stage cost $L\,:\, \mathbb S\,\times\, \mathbb A \rightarrow \mathbb R$, the closed-loop performance of a policy $\vect \pi\,:\, \mathbb S\,\rightarrow\, \mathbb A$ is given by:
\begin{align}
\label{eq:MDPCost}
J\left(\vect\pi\right) = \mathbb{E}_{\rho^{\vect \pi}}\left[\left.\sum_{k=0}^\infty \gamma^k L\left(\vect{s}_k,\vect a_k\right)\,\right |\, \vect a_k =\vect\pi\left(\vect{s}_k\right) \right],
\end{align}
where $\gamma\in (0,1]$ is a discount factor modeling the probabilistic lifetime of the system and ensuring the boundedness of the infinite-horizon problem \eqref{eq:MDPCost}. The expectation $\mathbb{E}_{\rho^{\vect \pi}}$ is taken over the Markov Chain resulting from~\eqref{eq:system} in closed-loop with policy $\vect\pi$ and initial conditions $\vect s_0\sim \rho_0$. The solution to the MDP provides an optimal policy $\vect\pi^\star$ from the set $\Pi$ of all admissible policies by minimizing \eqref{eq:MDPCost}, i.e.,
\begin{align}
\label{eq:OptPolicy}
\vect\pi^\star \in \argmin_{\vect\pi\in\Pi}\, J\left(\vect\pi\right).
\end{align}
% In this paper, we assume that \eqref{eq:system}, $L$ and $\Pi$ are such that \eqref{eq:MDPCost} is integrable. 
% \arash{Not very clear. Do you mean:
% \begin{align}
%     \Pi:=\{\vect\pi\,|\,  J\left(\vect\pi\right)<\infty\}.
% \end{align}
% I really do not think we need to define such a set for policies. The policies that yield infinite performance, clearly are not optimal. 
% }
The solution to the MDP is characterized by the Bellman equations: 
\begin{subequations}
\label{eq:Bellman}
\begin{align}
Q^\star\left(\vect s,\vect a\right) &= L\left(\vect s,\vect a\right) + \gamma \mathbb{E}_{\rho}\left[ V^\star\left(\vect{s}_+\right)\,|\, \vect s,\vect a \right] \, ,\\
V^\star\left(\vect s\right) &= \min_{\vect a}\,\, Q^\star\left(\vect s,\vect a\right) \label{eq:Vstar:Bell} \, , \\
\vect\pi^\star\left(\vect s\right) &\in \mathrm{arg} \min_{\vect a}\,\, Q^\star\left(\vect s,\vect a\right) \,, \label{eq:Pistar:Bell}
\end{align}
\end{subequations}
where $V^{\star}\, \text{and } Q^{\star}$ represent optimal value and optimal action-value functions respectively, and $\mathbb{E}_{\rho}\left[\cdot \mid \vect{s}, \vect{a}\right]$ refers to the expectation taken over the true system distribution \eqref{eq:system}. In the MDP framework, the constraints on the state and action can be enforced by assigning infinite values to the stage cost $L$ for constraint violation. More specifically, if a stage cost $L\left(\vect s,\vect a\right)$ is to be minimized by the optimal policy $\vect\pi^\star$ while respecting a set of constraints $\vect h\left(\vect s,\vect a\right)\leq 0$, a modified stage cost of the form:
\begin{align}
 \ell\left(\vect s,\vect a\right) = L\left(\vect s,\vect a\right) + \cmatr{ccc}{0&\mathrm{if}&\vect h\left(\vect s,\vect a\right) \leq 0\\
\infty&\mathrm{if}&\vect h\left(\vect s,\vect a\right) > 0 }, 
\end{align}
is used. We finally observe that equations \eqref{eq:Bellman} are generally very difficult to solve due to the curse of dimensionality of dynamic programming.

Model-based optimization frameworks such as MPC aiming to find an optimal policy for the true system \eqref{eq:system} typically consider a predictive model of the form:
\begin{align}
\label{eq:model}
\hat{\vect{s}}_+ \sim \hat{\rho} \left(\,.\,|\,\vect s, \vect a\,\right)\,,
\end{align}
as an approximation to the true system dynamics \eqref{eq:system}. Here $\hat{\rho}$ is the state transition probability under the model. Let $\hat{\vect \pi}^\star, \hat{Q}^\star, \hat{V}^\star$ be the optimal policy, optimal action-value function, and optimal value function, respectively, resulting as a solution to solving the MDP \eqref{eq:Bellman} under the model's Markov chain \eqref{eq:model}.  This approximate solution of the MDP based on the model (model-based MDP) typically does not match the optimal solution $\vect {\pi^\star}$ of the MDP \eqref{eq:OptPolicy}, derived based on the true system dynamics \eqref{eq:system}. In the remainder of this paper, we will establish the necessary and sufficient conditions on the model \eqref{eq:model} s.t. $\hat{\vect \pi}^\star = \vect {\pi^\star}$. MPC is a widely studied and significant special case of MDP. In the following, we provide the necessary background to establish the necessary and sufficient conditions for MPC to achieve an optimal policy. For this analysis, we consider MPC as a special case of an MDP with a deterministic predictive model.

\subsection{Model Predictive Control}\label{sec:MPC}

MPC approximately solves the optimization problem defined in \eqref{eq:MDPCost} using the model \eqref{eq:model}, resulting in a policy $\hat{\vect{\pi}}^\star$ which ideally should match the optimal policy $\vect{\pi^\star}$. Though any MPC formulation can, in principle, be used, we will focus on deterministic MPC in this paper. A deterministic MPC scheme solves the following optimization problem, for any given state $\vect{s}_k$ at time $k$:
\begin{subequations}
\label{eq:MPC0}
\begin{align}
V^\mathrm{MPC}\left(\vect{s}_k\right)=\min_{\hat{\vect{s}},\vect {u}}&\quad \gamma^{N} T\left(\hat{\vect{s}_{N}}\right) + \sum_{i=0}^{N-1} \gamma^k L\left(\hat{\vect{s}}_i,\vect{u}_i\right)\label{eq:MPC0:Cost}\\
\mathrm{s.t.} &\quad \hat{\vect{s}}_{i+1} = \vect f\left(\hat{\vect{s}}_i,\vect {u}_i\right)\label{eq:MPC0:Dyn} \,,\\
&\quad \vect h\left(\hat{\vect{s}}_i,\vect {u}_i\right)\leq 0 \label{eq:MPC0:Const} \,, \\
&\quad \hat{\vect{s}}_0 = \vect{s}_k,\quad \hat{\vect{s}}_N \in \mathbb T \,. \label{eq:MPC0:Boundaries} 
\end{align}
\end{subequations}
 Here $\vect f$ is the deterministic prediction model, $\vect {u}$ is the MPC actions/inputs and $\vect h$ represents a set of constraints to be respected. The terminal set $\mathbb T$ and the terminal cost $T$ are typically used to account for the fact that Problem~\eqref{eq:MPC0} is solved over a finite horizon $N$, while the actual MDP \eqref{eq:MDPCost} is open-ended. We define,  $V^\mathrm{MPC} = \infty$ when \eqref{eq:MPC0} is infeasible. We will follow the MPC formulation with discounting \eqref{eq:MPC0}, but all the results presented in this paper are valid also for undiscounted MPC. Note that a strong theory supporting the use of undiscounted MPC to support the functions related to discounted MDPs is provided in~\cite{Zanon2022b,Kordabad2023}.
 
 The deterministic MPC prediction model $\vect f$ is a special form of \eqref{eq:model},  which can be in a classic state-space form or in a data-driven prediction form. A typical type of model used in MPC is a simple expectation over the true system state transitions as follows:
\begin{align}
\label{eq:E:Fitting}
\hat{\vect{s}}_+ = \vect f\left(\vect s,\vect a\right) = \mathbb{E}_{\rho}\left[\vect{s}_+\,|\,\vect s,\vect a \right],\quad \forall\, \vect s,\vect a \,.
\end{align}
As we will discuss in the remainder of this paper, this is not necessarily the best model choice. Note that this is a specific form of the model from the general form in \eqref{eq:model}, other forms exist, such as those based on Maximum Likelihood Estimates (MLE), Bayesian estimates, etc. Stochastic versions of MPC use a stochastic model, which can be written as $\hat{\vect{s}}_+ = \vect f(\vect s, \vect a, \vect{\epsilon})$, where epsilon represents the stochasticity.

The solution of \eqref{eq:MPC0} delivers an input sequence $\vect {u}^\star_{0,\ldots,N-1}$ and the corresponding state sequence $\hat{\vect{s}}^\star_{0,\ldots,N}$ that the system is predicted to follow, according to the MPC model $\vect f$. Because the model is typically inexact, this prediction is usually inexact. To address this, MPC \eqref{eq:MPC0} is solved at every discrete time $k$, using the latest state of the system $\vect{s}_k$. Importantly, only the first input of the sequence, i.e., $\vect {u}^\star_{0}$ is applied to the true system. At the next time step, this process is repeated by updating the initial state of the MPC model \(\vect{f}\) with the true system state \(\vect{s}_{k+1}\), re-solving the MPC problem \eqref{eq:MPC0}, and applying only the first input from the newly computed sequence to the true system. Therefore, the MPC scheme  \eqref{eq:MPC0} produces the deterministic policy 
\begin{align}
\label{eq:MPC:Policy}
\vect\pi^\mathrm{MPC}\left(\vect{s}_k\right) = \vect {u}^\star_0,
\end{align}
which assigns to every feasible state $\vect{s}_k$ a corresponding action $\vect {u}^\star_0$ implicitly defined by \eqref{eq:MPC0}. The action-value function of the MPC is denoted by $Q^\mathrm{MPC}\left(\vect s,\vect a\right)$, is defined as:
\begin{subequations}
\label{eq:MPC:Qmodel}
\begin{align}
Q^\mathrm{MPC}\left(\vect{s}_k,\vect a_k\right):=\min_{\hat{\vect{s}},\vect a}&\quad 
\eqref{eq:MPC0:Cost} \label{eq:MPC:Qmodel:Cost}\\
%\gamma^{N} T\left(\vect x_{N}\right) + \sum_{i=0}^{N-1} \gamma^i L\left(\vect x_i,\vect {u}_i\right) \label{eq:MPC:Qmodel:Cost}\\
\mathrm{s.t.} &\quad \eqref{eq:MPC0:Dyn}-\eqref{eq:MPC0:Boundaries} \label{eq:MPC:Qmodel:const}\,,\\ 
%\vect x_{i+1} = \vect f\left(\vect x_i,\vect {u}_i\right)\label{eq:MPC0:Dyn} \, ,\\
%&\quad \vect h\left(\vect x_i,\vect {u}_i\right)\leq 0 \label{eq:MPC0:Const} \,, \\
%&\quad \vect x_{0} = \vect{s}_k,\quad \vect x_{N} \in \mathbb T \,, \\
&\quad  \vect {u}_0=\vect a_k. \label{eq:MPC:Qmodel:Aconst}
\end{align}
\end{subequations}
This definition of $Q^\mathrm{MPC}$ is valid in the sense of fundamental Bellman relationships between optimal action-value functions, value functions, and policies, i.e. 
\begin{subequations}
\begin{align}
V^\mathrm{MPC}\left(\vect s\right) &= \phantom{\mathrm{arg}}\min_{\vect a}\, Q^\mathrm{MPC}\left(\vect s,\vect a\right),\\\quad \vect \pi^\mathrm{MPC}\left(\vect s\right) &\in \mathrm{arg}\min_{\vect a}\, Q^\mathrm{MPC}\left(\vect s,\vect a\right),
\end{align}
\end{subequations}
hold by construction. Note that it has been established in~\cite{reinhardt2024economic} that an MPC scheme can be seen as a special case of MDP. We will exploit this fact in our proofs below.

\section{Optimality Conditions for MPC}\label{sec:conditons}
To establish the necessary and sufficient conditions on the predictive model for an MPC scheme to achieve optimality, we will first frame the results in the broader context of MDPs. We then translate these conditions to the MPC framework, given that MPC can be viewed as a degenerate form of an MDP. Throughout this note, we refer to the MDP under model \eqref{eq:model} as the ``model-based MDP'' and the MDP under the true (or real) system as the ``true MDP''. 

Consider the Bellman equation of the model-based MDP is given by:
\begin{subequations}
\begin{align}
    \hat{Q}^\star (\vect s, \vect a) &= L (\vect s, \vect a) + \gamma \mathbb{E}_{\hat{\rho}}\left[\hat{V}^{\star}\left(\vect{\hat{s}}_{+}\right) \mid \vect{s}, \vect{a}\right] \, ,\\
    \hat{\vect\pi}^\star\left(\vect s\right) &\in \mathrm{arg} \min_{\vect a}\,\, \hat{Q}^\star\left(\vect s,\vect a\right) \,, 
\end{align}
\end{subequations}
where \(\hat Q^\star\) and \(\hat V^\star\) are the optimal action-value and value functions respectively,
 $\vect{\hat{\pi}}^\star$ is the associated optimal policy, and $\mathbb{E}_{\hat{\rho}}\left[\cdot \mid \vect{s}, \vect{a}\right]$ denotes the expectation over the model distribution \eqref{eq:model}. We now introduce the following technical assumption which is central to the optimality conditions presented in this section.

\begin{Assumption}
\label{assum:technical_assumption}
The set:
\begin{equation} \label{eq:bounded_V}
\Omega := \Bigg\{ \vect s \in \mathcal{S} \, \Bigg |\,\left|\mathbb{E}_{\hat{\rho}}\left[\hat{V}^{\star}\left(\vect{\hat{s}}_k^{\pi^{\star}}\right)\right]\right|<\infty,  \quad \forall\, k < N \Bigg\}
\end{equation}
is assumed to be non-empty for stochastic trajectories of the model \((\vect{\hat{s}}_0^{\pi^{\star}}, \dots, _N)\) under the optimal policy \(\vect\pi^\star\).
\end{Assumption}

This assumption requires the existence of a non-empty set such that the optimal value function $\hat{V}^\star$ of the predicted optimal trajectories $(\vect{\hat{s}}_0^{\pi^{\star}},...,_{N})$ on the system model is finite for all initial states starting from this set. Assumption \ref{assum:technical_assumption} ensures that the predictive model trajectories under the optimal policy ${\vect\pi}^\star$ are contained within the set $\mathcal{S}$ where the value function ${V}^\star$ is bounded. This can be interpreted as a form of stability condition on $\hat{\rho}$ under the optimal trajectory $(\vect{\hat{s}}_0^{\pi^{\star}},...,_{N})$ \cite{gros2019data}.

We now introduce conditions for equivalence of the solution to the MDP with a model of the dynamics \eqref{eq:model} and the solution to the true MDP \eqref{eq:system} within the Bellman optimality framework \eqref{eq:Bellman} \cite{anand2024data}. This condition can be represented in terms of optimal action-value function as follows. If
\begin{align}
\label{eq:DT:PerfectModel}
\hat Q^\star\left(\vect s,\vect a\right)=Q^\star\left(\vect s,\vect a\right),\quad \forall\,\vect s,\vect a,
\end{align}
and Assumption~\ref{assum:technical_assumption} holds, then the solution to the MDP based on the model will yield the optimal policy, value function, and action-value function for the true MDP. In this sense, the MDP based on the model fully represents the true MDP. However, since we are interested in the optimality of the resulting policy, the condition \eqref{eq:DT:PerfectModel} can be relaxed such that only policy $\hat{\vect \pi}^\star$ needs to match the optimal policy $\vect\pi^\star$, not necessarily the value and action-value function. Then $\hat{Q}^\star$ needs to match $Q^{\star}$ only in the sense of:
 \begin{align}
\label{eq:DT:PerfectModel:argmin}
\mathrm{arg}\min_{\vect a}\,\hat{Q}^\star\left(\vect s,\vect a\right)=\mathrm{arg}\min_{\vect a}\,Q^\star\left(\vect s,\vect a\right),\quad \forall\,\vect s.
\end{align}
Note that \eqref{eq:DT:PerfectModel} implies \eqref{eq:DT:PerfectModel:argmin}, but the converse is not true, hence making \eqref{eq:DT:PerfectModel:argmin} less restrictive than \eqref{eq:DT:PerfectModel}. While less restrictive,  \eqref{eq:DT:PerfectModel:argmin} falls short of making the solution to the model-based MDP a fully satisfactory solution of the true MDP \eqref{eq:OptPolicy} in the sense of $V^\star$ and $Q^\star$.

\subsection{Necessary and Sufficient Condition} \label{sec:theorem}
Consider the following modification to $\hat{V}^\star$ of the model-based MDP, 
\begin{equation}
    {\hat{V}_{\lambda}}^\star(\vect s) = \hat{V}^\star(\vect s) + \lambda(\vect s) \,,
\end{equation}
which results in,
\begin{equation}\label{eq:modified_Q}
    {\hat{Q}_{\lambda}}^\star(\vect{s}, \vect{a}) = \hat{Q}^\star(\vect{s}, \vect{a})  + \lambda(\vect{s}) \,.
\end{equation}
Here, $\lambda(\vect{s})$ is a bounded function that can be chosen arbitrarily without affecting the optimal policy $\hat{\vect \pi}^\star$ of the model-based MDP, i.e.,
\begin{equation}
    \hat{\vect \pi}^\star \in \argmin_{\vect a} {\hat{Q}_{\lambda}}^\star(\vect{s}, \vect{a})  = \argmin_{\vect a} \hat{Q}^\star(\vect{s}, \vect{a}) + \lambda(\vect s)\,.
\end{equation} 
Defining,
\begin{equation}\label{eq:Lambda}
    \Lambda (\vect s, \vect a) := \lambda(\vect s) - \gamma \mathbb{E}_{\rho} [\lambda(\hat{\vect{s}}_+) \mid \vect{s}, \vect{a}]\,,
\end{equation}
the modified model-based MDP satisfies the Bellman equality:
\begin{equation}\label{eq:modified_bellman}
    {\hat{Q}_{\lambda}}^\star(\vect{s}, \vect{a}) = L(\vect{s}, \vect{a}) + \Lambda(\vect{s}, \vect{a}) + \gamma \mathbb{E}_{\hat{\rho}} [{\hat{V}_{\lambda}}^\star(\hat{\vect{s}}_+ \mid \vect{s}, \vect{a})] \,,
\end{equation}
given that Assumption~\ref{assum:technical_assumption} holds. In addition, we observe that the choice of $\lambda$ does not affect the optimal advantage function of the model-based MDP defined as
\begin{equation}
    {\hat{A}_{\lambda}}^\star(\vect s) := {\hat{Q}_{\lambda}}^\star(\vect s) - {\hat{V}_{\lambda}}^\star(\vect s).
\end{equation}

% The necessary optimality condition, if it exists,{\color{blue}\textbf{WHAT DOES THIS MEAN?? CONDITION OR SOLUTION? I also don't get the comment on (15)}} should match the optimal policy under true and model-based MDPs satisfying \eqref{eq:DT:PerfectModel:argmin}. 

In the following, for the sake of simplicity, let us assume that $Q^\star$ and ${\hat{Q}_{\lambda}}^\star$ are bounded on a compact set and infinite outside of it. 
Consider the advantage functions associated with the true MDP:
\begin{align}\label{eq:advantage_functions}
    A^{\star}(\mathbf{s}, \mathbf{a}) =Q^{\star}(\mathbf{s}, \mathbf{a})-V^{\star}(\mathbf{s}) \, .
\end{align}
By construction, it holds that
\begin{equation}
    \min _{\mathbf{a}} {\hat{A}_{\lambda}}^\star(\mathbf{s}, \mathbf{a})=\min _{\mathbf{a}} A^{\star}(\mathbf{s}, \mathbf{a})=0 \quad \forall \vect s  \,.
\end{equation}
With these preliminaries in place, the following two lemmas will help us characterize the optimality condition.

\begin{Lemma}\label{lem:necessary_condition_1}
    The existence a class $\mathcal{K}$ function~\footnote{A function $\alpha\,:\,\mathbb R_{\geq 0}\mapsto \mathbb R_{\geq 0}$ is said to belong to class $\mathcal{K}$ function if $\alpha(0)=0$ and $\alpha$ is a strictly increasing function.} $\alpha\,:\,\mathbb R_{\geq 0}\mapsto \mathbb R_{\geq 0}$ such that
    \begin{equation}\label{eq:aplha_for_A}
       {\hat{A}_{\lambda}}^\star(\mathbf{s}, \mathbf{a}) \geq \alpha (A^{\star}(\mathbf{s}, \mathbf{a})), \quad \forall\, \vect s\,, \vect a,
    \end{equation}
    is a necessary and sufficient condition for 
    \begin{equation}\label{eq:A_min_subset}
        \underset{\vect a}{\arg \min } {\hat{A}_{\lambda}}^\star(\mathbf{s}, \mathbf{a}) \subseteq \underset{\vect a}{\arg \min } A^{\star}(\mathbf{s}, \mathbf{a}), \quad \forall \vect s
    \end{equation}
    to hold.
\end{Lemma}
\begin{pf}
    We will prove the claim by  establishing that: (i) $\neg$ \eqref{eq:A_min_subset} $\Rightarrow$ $\neg$ \eqref{eq:aplha_for_A}; and (ii) \eqref{eq:A_min_subset} $\Rightarrow$ \eqref{eq:aplha_for_A}. 
    
    Condition $\neg$ \eqref{eq:A_min_subset}  implies that there exist a state-action pair $\overline{\vect s}, \overline{\vect a} $ such that
    \begin{equation}
        \bar{\vect a}\in \argmin_{\vect a}  {\hat{A}_{\lambda}}^\star\left(\bar{\vect s},\vect a\right) ,\quad\text{and}\quad \bar{\vect a}\notin \argmin_{\vect a} A^\star\left(\bar{\vect s},{\vect a}\right),
        \end{equation}
i.e.,        
    \begin{equation}
        {\hat{A}_{\lambda}}^\star(\overline{\mathbf{s}}, \overline{\mathbf{a}})=0, \quad \text { and } \quad A^{\star}(\overline{\mathbf{s}}, \overline{\mathbf{a}})>0 \,.
    \end{equation}
    Then Condition \eqref{eq:aplha_for_A} can not hold for any $\mathcal{K}$ function $\alpha$ at $\overline{\vect s}, \overline{\vect a}$. 

    Next we establish that \eqref{eq:A_min_subset} $\Rightarrow$ \eqref{eq:aplha_for_A}. To that end, we first construct the following function $\alpha_0$: 
    \begin{subequations}
            \label{eq:alpha:K}
            \begin{align}
                \alpha_0(x) :=
                \min_{ {\vect s, \vect a}} \quad &{\hat{A}_{\lambda}}^\star\left({\vect s}, {\vect a}\right)  \\
                \mathrm{s.t.}\quad & A^\star\left({\vect s}, {\vect a}\right) \geq  x\, ,
            \end{align}
        \end{subequations}
     for all $x\leq \bar A^\star$, where 
     \begin{align}
        \bar A^\star=\max_{ {\vect s, \vect a}} A^\star (\vect s,\vect a)\,,
     \end{align}
    and
     \begin{align}\label{eq:alpha:K:2}
         \alpha_0(x) := \alpha_0(\bar A^\star)+x-\bar A^\star  \quad \forall\, x>\bar A^\star\,.
     \end{align}
    We now prove that $\alpha_0$ satisfies the following properties:
    \begin{enumerate}
       \item\label{item:alpha0:1} $\alpha_0(0)=0$.
        \item \label{item:alpha0:2} $\alpha_0(x)>0$ for $x>0$.
        \item \label{item:alpha0:3} $\alpha_0(x)$ is a (not necessarily strictly) increasing function.
        \item \label{item:alpha0:4} $\alpha_0(x)$ satisfies \eqref{eq:aplha_for_A}.
    \end{enumerate}
    For condition \ref{item:alpha0:1}), the feasible domain of~\eqref{eq:alpha:K} for $x=0$, contains all $\vect s, \vect a$, including $(\vect s,\hat{\vect\pi}^\star (\vect s))$, therefore  $\alpha_0(0)=0$. 
    
    We prove condition \ref{item:alpha0:2}) using contradiction. First, we observe that $\alpha_0(x)$ is non-negative for all $x$. Assume that $\alpha_0 (x)=0$ for some $x\in (0, \bar A^\star]$, then there exists $\bar{\vect s}$ and $\bar{\vect a}$ such that ${\hat{A}_{\lambda}}^\star\left(\bar {\vect s}, \bar {\vect a}\right)=0$, while respecting constraint  $A^\star\left(\bar {\vect s}, \bar {\vect a}\right)\geq x>0$. However, this contradicts \eqref{eq:A_min_subset}. Additionally, for all $x>\bar A^\star$, $\alpha_0(x)> 0$ by construction, see~\eqref{eq:alpha:K:2}.
    
    To prove \ref{item:alpha0:3}), we first define the set of feasible solutions of~\eqref{eq:alpha:K} at a given $x$ by $\mathfrak C (x)$. One can observe that $\mathfrak C (x_2)\subseteq \mathfrak C (x_1)$ for any $x_1 \leq x_2$, making the optimal objective (or $\alpha_0(x)$) non-decreasing when $x$ increases. For all $x>\bar A^\star$, we have $\alpha_0(x)\geq \alpha_0(\bar A^\star)$, and $\alpha_0(x)$ is also increasing by construction in~\eqref{eq:alpha:K:2}.
    
    We then prove \ref{item:alpha0:4}), i.e., we show that, for all $\vect s_0, \vect a_0$, we have:
     \begin{align}
     \label{eq:Algebraicy:K0}
     \hat A^\star\left(\vect s_0,\vect a_0\right) \geq \alpha_0\left( A^\star\left(\vect s_0,\vect a_0\right)\right).
    \end{align}
    Since $A^\star\left(\vect s_0,\vect a_0\right)\leq \bar A^\star$,  from~\eqref{eq:alpha:K}, we have:
    \begin{subequations}
            \label{eq:alpha:K0}
            \begin{align}
                \alpha_0(A^\star\left(\vect s_0,\vect a_0\right)) :=
                \min_{{\vect s}, {\vect a}} \quad &\hat A^\star\left({\vect s}, {\vect a}\right)  \\
                \mathrm{s.t.} \quad & A^\star\left({\vect s}, {\vect a}\right) \geq  A^\star\left(\vect s_0,\vect a_0\right).
            \end{align}
        \end{subequations}
    Note that $\vect s=\vect s_0$, $\vect a=\vect a_0$ is a feasible pair for~\eqref{eq:alpha:K0}, such that optimality directly yields~\eqref{eq:Algebraicy:K0}. Furthermore, one can verify that any function satisfying conditions~\ref{item:alpha0:1}), \ref{item:alpha0:2}) and \ref{item:alpha0:3}) can be lower bounded by a $\mathcal{K}$ function, concluding the proof.
    \end{pf}

\begin{Lemma} \label{lem:necessary_condition_2} The existence of a $\mathcal{K}$ function $\beta\,:\,\mathbb R_{\geq 0}\mapsto \mathbb R_{\geq 0}$ such that
    \begin{align}
    \label{eq:Algebraicy:rev:K}
    \beta(A^\star\left(\vect s,\vect a\right)) \geq {\hat{A}_{\lambda}}^\star\left(\vect s,\vect a\right) ,\quad \forall\, \vect s,\vect a
   \end{align}
   is a necessary and sufficient condition for
   \begin{align}
   \label{eq:AMinSetMatchy:rev:K}
   \argmin_{\vect a} A^\star\left(\vect s,\vect a\right) \subseteq \argmin_{\vect a}{\hat{A}_{\lambda}}^\star \left(\vect s,\vect a\right),\quad \forall\, \vect s
   \end{align}
   to hold.
    \end{Lemma}

   \begin{pf} 
   We first prove that the claim of the lemma is equivalent to the following. The existence of a $\mathcal{K}$ function $\tilde\beta\,:\,\mathbb R_{\geq 0}\mapsto \mathbb R_{\geq 0}$ such that
    \begin{align}
    \label{eq:Algebraicy:K2}
    A^\star\left(\vect s,\vect a\right) \geq \tilde\beta\left( {\hat{A}_{\lambda}}^\star \left(\vect s,\vect a\right)\right),\quad \forall\, \vect s,\vect a
   \end{align}
   is a necessary and sufficient condition for
   \begin{align}
   \label{eq:AMinSetMatchy:K2}
   \argmin_{\vect a} A^\star\left(\vect s,\vect a\right) \subseteq \argmin_{\vect a} {\hat{A}_{\lambda}}^\star \left(\vect s,\vect a\right),\quad \forall\, \vect s
   \end{align}
   to hold. 
   
   Since $\tilde\beta$ is a $\mathcal{K}$ function, then there exists a $\mathcal{K}$ function $\beta\,:\,\mathbb R\to\mathbb R$ satisfying~\cite{sontag1989smooth}
   % \mario{ maybe we should cite something here}
   \begin{align}
   \beta\left(\tilde \beta\left(x\right)\right) = x.
    \end{align}
   Then applying $\beta$ to both sides of the inequality  \eqref{eq:Algebraicy:K2} implies that  \eqref{eq:Algebraicy:K2} is equivalent to \eqref{eq:Algebraicy:rev:K}.
   
Now, the proof of Lemma~\ref{lem:necessary_condition_2} directly follows from the proof of Lemma~\ref{lem:necessary_condition_1}, utilizing the equivalence between conditions \eqref{eq:Algebraicy:K2}–\eqref{eq:AMinSetMatchy:K2} and \eqref{eq:aplha_for_A}–\eqref{eq:A_min_subset} in Lemma~\ref{lem:necessary_condition_1}.

   \end{pf}
Additionally, the following Corollary follows directly from Lemma~\ref{lem:necessary_condition_1} and Lemma~\ref{lem:necessary_condition_2}.

\begin{corollary}
    The existence of class $\mathcal{K}$ functions $\alpha, \beta$ s.t.
    \begin{equation}\label{eq:necesary_condition_A}
        \beta (A^{\star}(\mathbf{s}, \mathbf{a})) \geq {\hat{A}_{\lambda}}^\star(\mathbf{s}, \mathbf{a}) \geq \alpha (A^{\star}(\mathbf{s}, \mathbf{a})), \,\,\forall \vect s,\vect a
    \end{equation}
    is a necessary and sufficient condition for
    \begin{equation}\label{eq:necesary_condition_Amin}
        \underset{{\vect a}}{\arg \min } {\hat{A}_{\lambda}}^\star(\mathbf{s}, \mathbf{a})=\underset{{\vect a}}{\arg \min } A^{\star}(\mathbf{s}, \mathbf{a}), \quad \forall \mathbf{s}\,
    \end{equation}
    to hold.
\end{corollary}
\begin{pf}

The proof is derived by combining the conditions \eqref{eq:Algebraicy:rev:K} and \eqref{eq:AMinSetMatchy:rev:K} from Lemma~\ref{lem:necessary_condition_2} with the conditions \eqref{eq:aplha_for_A} and \eqref{eq:A_min_subset} from Lemma~\ref{lem:necessary_condition_1}, resulting in the left and right sides of the inequality \eqref{eq:necesary_condition_A} and the condition \eqref{eq:necesary_condition_Amin}.
    
\end{pf}

\begin{thm}
    \label{Th:necesary_condition}
    The necessary and sufficient condition for a model-based MDP, defined under the model state transition distribution \eqref{eq:model}, to yield the optimal policy $\vect \pi^\star$ of the true MDP defined under the true system state transition distribution \eqref{eq:system}, is as follows:
    \begin{equation}\label{eq:necessary_condition}
        \begin{aligned}
        \beta (L & (\mathbf{s}, \mathbf{a})+ \gamma \mathbb{E}_{\rho}\left[V^{\star}\left(\mathbf{s}_{+}\right) \mid \mathbf{s}, \mathbf{a}\right]- V^{\star}(\mathbf{s})) \geq \\
        & L(\mathbf{s}, \mathbf{a})+ \Lambda(\mathbf{s}, \mathbf{a})+\gamma \mathbb{E}_{\hat{\rho}}\left[{V}^{\star}\left(\mathbf{\hat{s}}_{+}\right) \mid \mathbf{s}, \mathbf{a}\right]-{V}^\star(\mathbf{s}) \geq \\
        & \hspace{1cm} \alpha (L(\mathbf{s}, \mathbf{a})+\gamma \mathbb{E}_{\rho}\left[ V^{\star}\left(\mathbf{s}_{+}\right) \mid \mathbf{s}, \mathbf{a}\right]- V^{\star}(\mathbf{s}))
        \end{aligned}
    \end{equation}
    for all $\vect s, \vect a$ with
    \begin{align}
    \Lambda\left(\vect s,\vect a\right) := \lambda\left(\vect s\right) - \gamma  \mathbb{E}_{\rho}[\lambda(\hat{\vect{s}}_+ \mid \vect{s}, \vect{a})], \label{eq:storage}
    \end{align}
    and $\lambda\left(\vect s\right)$ chosen such that
    \begin{align}
    \label{eq:ValueFunctionMatching}
    {\hat{V}_{\lambda}}^\star\left(\vect s \right) = V^\star\left(\vect s \right)\,,
    \end{align} 
    where $V^\star$ and ${\hat{V}_{\lambda}}^\star$ are bounded on the set $\Omega$ as defined in Assumption \ref{assum:technical_assumption}. 
\end{thm}

\begin{pf}
    The necessary condition \eqref{eq:necesary_condition_A} entails
    \begin{align}
        \beta({Q}^\star(\mathbf{s}, \mathbf{a})- {V}^\star(\mathbf{s})) &\geq \, {\hat{Q}_{\lambda}}^\star(\mathbf{s}, \mathbf{a}) -{\hat{V}_{\lambda}}^\star(\mathbf{s}) \notag \\
        &\geq \alpha (Q^{\star}(\mathbf{s}, \mathbf{a})- V^{\star}(\mathbf{s}))\,. 
    \end{align}
    Hence, using \eqref{eq:Bellman} and \eqref{eq:modified_bellman}, we arrive at the following inequality:
    \begin{equation}\label{eq:condition_without_V_star}
    \begin{aligned}
    \beta (L & (\mathbf{s}, \mathbf{a})+ \gamma \mathbb{E}_{\rho}\left[V^{\star}\left(\mathbf{s}_{+}\right) \mid \mathbf{s}, \mathbf{a}\right]- V^{\star}(\mathbf{s})) \geq \\
    & L(\mathbf{s}, \mathbf{a})+ \Lambda(\mathbf{s}, \mathbf{a})+\gamma \mathbb{E}_{\hat{\rho}}\left[ {\hat{V}_{\lambda}}^{\star}\left(\mathbf{\hat{s}}_{+}\right) \mid \mathbf{s}, \mathbf{a}\right]-{\hat{V}_{\lambda}}^\star(\mathbf{s}) \geq \\
    & \hspace{1cm} \alpha (L(\mathbf{s}, \mathbf{a})+\gamma \mathbb{E}_{\rho}\left[ V^{\star}\left(\mathbf{s}_{+}\right) \mid \mathbf{s}, \mathbf{a}\right]- V^{\star}(\mathbf{s}))\,.
    \end{aligned}
    \end{equation}
    Since ${\hat{A}_{\lambda}}^\star$ is independent of $\lambda$, that function can be chosen arbitrarily without affecting the necessary condition \eqref{eq:necesary_condition_A} nor the functions $\alpha$ and $\beta$. It follows that (i) choosing a specific $\lambda$ does not weaken the necessary condition to a sufficient one, and (ii) the value of $\lambda$ can always be selected such that
    \begin{equation}\label{eq:lambda_for_alpha}
        {\hat{V}_{\lambda}}^\star(\vect s) = V^\star (\vect s)\,, \quad \forall \vect s\,,
    \end{equation}
    with ${\hat{V}_{\lambda}}^\star(\vect s)$  and $V^\star(\vect s)$ being bounded on the same set, leads to the necessary condition \eqref{eq:necessary_condition}. 
\end{pf}

\begin{corollary}[Sufficient Condition]
    \label{Th:sufficient_condition}
    Assume that $\alpha$ and $\beta$ are both the identity function and $\lambda = \Delta$ is constant for all $\vect s$, the necessary condition \eqref{eq:necessary_condition} reduces to a sufficient condition, given by:
    \begin{align}
    \label{eq:suficient_condition}
    \mathbb{E}_{\rho}\left[V^\star\left(\vect{s}_+\right)\,|\, \vect s,\vect a\, \right] - \mathbb{E}_{\hat{\rho}}\left[{V}^\star\left(\hat{\vect s}_+\right)\,|\, \vect s,\vect a\, \right]  = \Delta \,.
    \end{align}
\end{corollary}
\begin{pf} 
    Consider the optimal value function of the model-based MDP: 
    \begin{equation}
    {\hat{Q}_{\lambda}}^{\star} (\vect s, \vect a) = L(\vect s,\ \vect a) + \gamma \mathbb{E}_{\hat{\rho}}\left[{V}^\star\left(\hat{\vect s}_+\right)\,|\, \vect s,\vect a\, \right] \label{eq:V_pi_model}
    \end{equation}
    Using \eqref{eq:suficient_condition},  we get:
    \begin{align}
        {\hat{Q}_{\lambda}}^{\star} (\vect s, \vect a) &= L(\vect s,\ \vect a)  + \gamma \mathbb{E}_{\rho}\left[V^\star\left(\vect{s}_+\right)\,|\, \vect s,\vect a\, \right] - \gamma \Delta \\
        &={Q}^{\star} (\vect s, \vect a) - \gamma \Delta. \label{eq:DT:PerfectModel:PlusConstant}
    \end{align}
    % This can be rearranged as:
    % \begin{align}\label{eq:DT:PerfectModel:PlusConstant}
    %     \hat{Q}^{\star} (\vect s, \vect a) + \gamma \Delta = {Q}^{\star} (\vect s, \vect a) \,.
    % \end{align}
    Note that \eqref{eq:DT:PerfectModel:PlusConstant} can be converted to Condition \eqref{eq:DT:PerfectModel} simply by adding the constant $\gamma\Delta$ to the costs in \eqref{eq:MDPCost}, without affecting the optimal policy. Therefore, this condition is sufficient for the optimality of the model-based MDP.
\end{pf}

\begin{corollary}
    The necessary and sufficient condition on the deterministic predictive model $\hat \rho (\hat {\vect s}_+ | \vect s,\vect a) = \delta(\hat {\vect s}_+-\vect f(\vect s, \vect a))$, where $\delta$ is the Dirac delta, of an MPC scheme \eqref{eq:MPC0}, such that it provides the optimal policy $\vect \pi ^\star$ for the true MDP, is given by:
        \begin{equation}\label{eq:ConditionFull}
        \begin{aligned}
        \beta (L & (\mathbf{s}, \mathbf{a})+ \gamma \mathbb{E}_{\rho}\left[V^{\star}\left(\mathbf{s}_{+}\right) \mid \mathbf{s}, \mathbf{a}\right]- V^{\star}(\mathbf{s})) \geq \\
        & L(\mathbf{s}, \mathbf{a})+ \Lambda(\mathbf{s}, \mathbf{a})+\gamma {V}^\star\left(\vect f\left(\vect s,\vect a\right)\right) -{V}^\star(\mathbf{s}) \geq \\
        & \hspace{1cm} \alpha (L(\mathbf{s}, \mathbf{a})+\gamma \mathbb{E}_{\rho}\left[ V^{\star}\left(\mathbf{s}_{+}\right) \mid \mathbf{s}, \mathbf{a}\right]- V^{\star}(\mathbf{s}))
        \end{aligned}
    \end{equation}
    for all $\vect s, \vect a$ with $T = V^{\star}$. 
\end{corollary}
\begin{pf}
Consider the modified form of $Q^\mathrm{MPC}$ in \eqref{eq:MPC:Qmodel}, in line with the modification in \eqref{eq:modified_Q} with the function $\lambda (\vect s)$:
\begin{subequations}
\label{eq:MPCPolicy}
\begin{align}
{Q_{\lambda}}^\mathrm{MPC}\left(\vect s,\vect a\right) = \lambda\left(\vect s\right)+ \min_{\vect x,\vect {u}}&\,\, \gamma^NT\left(\vect x_N\right)+\sum_{k=0}^{N-1}\gamma^kL\left(\vect x_k,\vect {u}_k\right)\\
\mathrm{s.t.}&\quad \vect x_{k+1}=\vect f\left(\vect x_k,\vect {u}_k\right)\,, \\
&\quad \vect x_0=\vect s,\quad \vect {u}_0=\vect a\
\end{align}
\end{subequations}
MPC \eqref{eq:MPCPolicy} satisfies the modified Bellman equations equivalent to \eqref{eq:modified_bellman}:
\begin{subequations}\label{eq:modified_bellman_MPC}
\begin{align}
{Q_{\lambda}}^\mathrm{MPC}\left(\vect s,\vect a\right) &=  L\left(\vect s,\vect a\right) + \Lambda\left(\vect s,\vect a\right)+ \gamma {V_{\lambda}}^\mathrm{MPC}\left(\vect f\left(\vect s,\vect a\right)\right) \,.
\end{align}
\end{subequations}
The proof follows by combining 
    \begin{equation}\label{eq:value_model_to_value_mdp}
        \mathbb{E}_{\hat{\rho}}\left[{V}^\star\left(\hat{\vect{s}}_+\right)\,|\, \vect s,\vect a   \right] = {V}^\star\left(\vect f\left(\vect s,\vect a\right)\right)\,
    \end{equation}
and \eqref{eq:necessary_condition}. 
\end{pf}

\section{Discussions}\label{sec:discussions}
%on the condition
 We observe that the necessary and sufficient optimality conditions from Theorem \ref{Th:necesary_condition} interlace the optimal value function of the true MDP, $V^\star$, and the corresponding model-based MDP (or MPC) in a non-trivial way. Importantly, the necessary condition \eqref{eq:necessary_condition} or the sufficient condition \eqref{eq:suficient_condition} are unlike the conventional criteria for model estimation that best fit the data, for example as given in \eqref{eq:E:Fitting} or other estimation methods such as MLE or Bayesian estimation. Another key observation about the necessary and sufficient conditions in MPC is that they apply only to the single-step predictions of the model $\vect f$ in the MPC scheme, and not to the full prediction of the state trajectory inherent in MPC decisions. This is evident from the conditions in \eqref{eq:necessary_condition}, which are intrinsically linked to the single-step predictions of the MPC model. 
This suggests that the argument on optimality is relevant only to predictive models capable of making one-step-ahead predictions \cite{anand2024data}. 

%special case of expercation based models

It is important to note that among conventional data-fitted models, the expectation-based models in \eqref{eq:E:Fitting} can achieve local optimality under specific conditions as discussed in \cite{anand2024data}. One such class of problems includes ``economic'' MDPs with smooth dynamics and smooth objectives, achieving dissipativity \cite{gros2022economic} to a narrow attraction set. The second class comprises ``tracking'' problems with smooth dynamics of compact support, where the stage cost is quadratic \cite{anand2024data}. For cases outside these restricted classes of problems, an MPC scheme based on \eqref{eq:E:Fitting} can not guarantee even local optimality. Examples include MDPs that do not achieve dissipativity or those with non-smooth cost functions or dynamics \cite{anand2024data}. The necessary and sufficient conditions presented here clarify why an MPC scheme using a predictive model constructed through conventional data-fitting methods does not necessarily achieve optimality or improve performance in practice.

\section{Conclusion} \label{sec:conclusion}
We introduced the necessary and sufficient conditions for a predictive model such that the resulting MPC scheme achieves closed-loop optimality. These conditions are significantly different from the conventional approach of constructing MPC models for prediction accuracy, which typically focuses on the best fit to the data. These conditions underscore the importance of looking beyond the prediction accuracy of MPC models, particularly when achieving optimality is of importance.

\bibliographystyle{plain}        % Include this if you use bibtex 
\bibliography{main}           % and a bib file to produce the 
\end{document}